\documentclass{amsart}
\usepackage{amsfonts}
\usepackage{amsmath,amssymb}
\usepackage{amsthm}
\usepackage{amscd}
\usepackage{graphics}
\usepackage{graphicx}

\theoremstyle{remark}{
\newtheorem{Def}{{\rm Definition}}
\newtheorem{Ex}{{\rm Example}}
\newtheorem{Rem}{{\rm Remark}}
\newtheorem{Prob}{{\rm Problem}}

}
\theoremstyle{plain}
{

\newtheorem{MainThm}{Main Theorem}

}

\begin{document}
\title[Real algebraic functions whose Reeb spaces are given graphs]{Real algebraic functions on closed manifolds whose Reeb spaces are given graphs}
\author{Naoki kitazawa}
\keywords{Kronrod-Reeb graphs. Real algebraic manifolds. Real algebraic domains. Poincar\'e-Reeb graph. \\
\indent {\it \textup{2020} Mathematics Subject Classification}: Primary~57R45, 58C05. Secondary~57R19.}
\address{Institute of Mathematics for Industry, Kyushu University, 744 Motooka, Nishi-ku Fukuoka 819-0395, Japan\\
 TEL (Office): +81-92-802-4402 \\
 FAX (Office): +81-92-802-4405 \\
}
\email{n-kitazawa@imi.kyushu-u.ac.jp}
\urladdr{https://naokikitazawa.github.io/NaokiKitazawa.html}
\maketitle
\begin{abstract}

In our paper, we construct a real-algebraic function on some closed manifold whose {\it Reeb} ({\it Kronrod-Reeb}) {\it graph} is a graph respecting some {\it algebraic domain}: a graph for this is called a {\it Poincar\'e-Reeb} graph.

The {\it Reeb graph} of a smooth function is defined as a natural graph which is the quotient space of the manifold of the domain under a natural equivalence relation for some wide and nice class of smooth functions. 
The vertex set is defined as the set of all connected components containing some singular points of the function: a {\it singular} point of a smooth function is a point where the differential vanishes.
Morse-Bott functions give very specific cases.
The relation is
 to contract each connected component of each preimage to a point. 

Sharko has asked a natural and important problem: can we construct a nice smooth function whose Reeb graph is a given graph? Explicit answers have been given first by Masumoto-Saeki in a generalized manner for closed surfaces. After that various answers have been presented by various researchers and most of them are essentially for functions on closed surfaces and Morse functions such that connected components of preimages containing no singular points are spheres. Recently the author has also considered questions and solved in cases the preimages are general manifolds. 

\end{abstract}
\section{Introduction}
\label{sec:1}

The {\it Reeb graph} ({\it Kronrod-Reeb graph})
 of a differentiable function $c:X \rightarrow Y$ is a graph whose underlying topology is 
the quotient space of the manifold of the domain defined in the following way: two points $x_1, x_2 \in X$ are equivalent if and only if they are in a same connected component of some preimage $c^{-1}(y)$.

Let $W_c$ denote this. We define the quotient map $q_c:X \rightarrow W_c$ and we can define a map $\bar{c}$ enjoying the relation $c=\bar{c} \circ q_c$ uniquely. 

This is regarded as a graph for some wide class of smooth functions. We explain about the structure of our Reeb graph for $W_c$.
\begin{Rem}

We call such graphs "Reeb graphs" in considerable situations simply where "kronrod-Reeb graphs" are also well-known of course.
\end{Rem}
Such topological and combinatorial objects have been fundamental and strong tools in understanding the manifolds. \cite{reeb} is a related pioneering paper.

We define fundamental terminologies, notions, and notation. 

For a topological space $X$ having the structure of some cell complex whose dimension is finite, its dimension is uniquely defined as the dimension of the cell complex. Let $\dim X$ denote this.

A polyhedron and a CW complex is of course of such a class and a topological manifold is of such a class, being known to have the structure of a CW complex of a finite dimension. 

A graph is a $1$-dimensional CW complex where the vertex set is the set of all $0$-dimensional cells and the edge set is the set of all $1$-dimensional cells. If the closure of an edge is homeomorphic to a circle, then it is called a {\it loop}. Hereafter, we do not consider such graphs and a graph is always a $1$-dimensional polyhedron.
An {\it isomorphism} from a graph $K_1$ onto $K_2$ means a piecewise smooth homeomorphism mapping the edge set and the vertex set of $K_1$ onto those of $K_2$.

For a differentiable map $c:X \rightarrow Y$,  a {\it singular} point $x \in X$ is a point where the differential is smaller than both $\dim X$ and $\dim Y$. The value at some singular point of $c$ is a {\it singular value} of $c$.
If a smooth function on a closed manifold has finitely many singular values, then the Reeb space of it is a graph where the vertex set consists of all points $p$ whose preimages ${q_c}^{-1}(p)$ contain some singular points of $c$. 
Morse-Bott functions and smooth functions of some considerably wide classes satisfy this. This is due to \cite{saeki2}.

\begin{Prob}
\label{prob:1}
For a graph, can we construct a smooth function on some closed manifold whose Reeb graph is isomorphic to this and which enjoys some nice (differential) topological properties and properties on singularity? We do not need to fix the manifold of the domain.
\end{Prob}
\cite{sharko} has asked this first. Smooth functions on closed surfaces have been explicitly constructed for graphs satisfying some nice conditions there. \cite{masumotosaeki} generalizes this for arbitrary graphs.
Later, \cite{martinezalfaromezasarmientooliveira, michalak} have given answers. These study cases for smooth functions on closed surfaces and Morse functions such that each connected component of each preimage containing no singular points is always a sphere for example. In \cite{kitazawa4, kitazawa5, kitazawa8, kitazawa9}, the author has studied cases where such connected components are general manifolds with mild conditions on singularities of the functions for example. \cite{saeki2} appears as a paper motivated by \cite{kitazawa4} and through related informal discussions by us. 
\begin{Prob}
\label{prob:2}
Can we construct a smooth map in Problem \ref{prob:1} as a morphism from a nicer or finer category. In other words, can we construct this as a real analytic one, and as a nicer one, real algebraic one, for example? 
\end{Prob}
The category we discuss has been always the smooth one. In this paper, we try Problem \ref{prob:2} and obtain related answers under some specific or explicit situations.

To present our main result, we need some terminologies, notions and notation from the theory of smooth manifolds and (real) algebraic manifolds for example.

${\mathbb{R}}^k$ is the $k$-dimensional Euclidean space and this is a simplest smooth manifold of dimension $k$ for an arbitrary integer $k>0$. This is also
 a simplest real algebraic manifold, which is also called the ($k$-dimensional) {\it real affine space}. 
It is also a Riemannian manifold equipped with the so-called Euclidean metric. For a point $x \in {\mathbb{R}}^k$, $||x|| \geq 0$ is the distance between $x$ and
 the original $0 \in {\mathbb{R}}^k$. 
$\mathbb{R}:={\mathbb{R}}^1$.
$S^k:=\{x \in {\mathbb{R}}^{k+1} \mid ||x||=1 \}$ is the $k$-dimensional unit sphere. This is a $k$-dimensional algebraic submanifold of ${\mathbb{R}}^{k+1}$ which is compact and has no boundary. It is also a smooth submanifold. It is connected for $k \geq 1$ and it is a discrete set with exactly two points for $k=0$.

An {\it algebraic domain} $D$ of ${\mathbb{R}}^k$ is some open subset there such that the boundary of the closure $\overline{D}$ consists of finitely many smooth algebraic submanifolds of dimension $k-1$ or smooth algebraic hypersurfaces with no boundaries. 

To simplify our arguments, let us assume the following where $l \geq 0$ is an integer.
\begin{itemize}
\item For each hypersurface $S_j$ in the family $\{S_j\}_{j=1}^l$, a real polynomial $f_{{\rm P}, S_j}$ is given so that the zero set and $S_j$ coincide and that the polynomial function $f_{{\rm P}, S_j}:S_j \rightarrow \mathbb{R}$ defined canonically has no singular points on $S_j$.
\item $D$ is assumed to be the intersection ${\bigcap}_{j=1}^l  \{x \in {\mathbb{R}}^k \mid  f_{{\rm P}, S_j}(x)>0 \}$.
\end{itemize}

For example, the interior ${\rm Int}\ D^k$ of $D^k:=\{x \in {\mathbb{R}}^{k} \mid ||x|| \leq 1 \}$ is a simplest example and $D^k$ is the $k$-dimensional unit disk. This is also a $k$-dimensional smooth, compact and connected submanifold.
Note that $||x||={\Sigma}_{j=1}^k {x_j}^2$ where $x:=(x_1,\cdots,x_k)$.

A {\it Poincar\'e-Reeb graph} is defined for a pair of an algebraic domain $D$ of the real affine space of dimension $k>1$ and a canonical projection ${\pi}_{k,1}$ mapping $(x_1,x_2) \in {\mathbb{R}}^{k}$ to $x_1 \in {\mathbb{R}}$. This can be presented in a more general manner.
Hereafter, we mainly respect the preprint \cite{bodinpopescupampusorea} and there such cases are discussed. Note that terminologies and situations are different in considerable cases and that here we can argue in a self-contained way.
\begin{Def}
\label{def:1}
A {\it Poincar\'e-Reeb graph for the pair $(D,{\pi}_{k,1})$} is a graph in the real affine space embedded by a piecewise smooth embedding with the following conditions.
\begin{itemize}
\item Each edge $e$ intersects each preimage of the projection ${\pi}_{k,1}$ in a so-called {\it generic} way or satisfying the "transversality". In other words, each edge is embedded smoothly and for each point $p_e$ in each edge $e$, the image of the differential at the point and the tangent vector space at the value $v(p_e)$ in the preimage ${{\pi}_{k,1}}^{-1}(p)$ of a suitable {\rm (}unique{\rm )} point $p$ by the projection ${\pi}_{k,1}$ generate the tangent vector space at the point $v(p_e) \in {\mathbb{R}}^k$. 
\item Two points in the closure $\overline{D}$ of $D$ can be defined to be equivalent if and only if they are in a same connected component of the preimage $\overline{D} \bigcap {{\pi}_{k,1}}^{-1}(p)$ for some point $p \in {\mathbb{R}}$ and the map obtained by the restriction of the projection to the closure $\overline{D}$. Let ${\pi}_{D}$ denote the restriction to  the closure $\overline{D}$.
Our Poincar\'e-Reeb graph for the pair can be also defined as the quotient space obtained by this equivalence relation. This is isomorphic to the Reeb graph of ${\pi}_{D}$. Furthermore, an isomorphism is defined as the canonically obtained correspondence. 
\item The vertex set of our Poincar\'e-Reeb graph for the pair is the union of the set of all singular points of the restrictions of the projection ${\pi}_{k,1}$ or ${\pi}_D$ to all connected components of the boundary $\partial \overline{D} \subset \overline{D}$. This set is also finite. 
\end{itemize}
\end{Def}
See also \cite{sorea1, sorea2} for related theory for example. We present our main result. In the following section we prove this and present related comments as our main content.
\begin{MainThm}
\label{mthm:1}
Consider a Poincar\'e-Reeb graph $K$ for the pair in Definition \ref{def:1} such that the closure $\overline{D}$ is compact. Take an arbitrary integer $k_0>k+1$. Then we can construct a real algebraic function on some {\rm (}$k_0-1${\rm )}-dimensional smooth closed manifold regarded as a smooth real algebraic manifold whose Reeb graph is isomorphic to the graph $K$ as a graph.
\end{MainThm}

\noindent {\bf Conflict of interest.} \\
	The author is a member of the project JSPS KAKENHI Grant Number JP22K18267 "Visualizing twists in data through monodromy" (Principal Investigator: Osamu Saeki). Our present study is supported by the project. \\
	\ \\
	{\bf Data availability.} \\
	Data essentially supporting our present study are all contained in our present paper.
\section{On Main Theorem \ref{mthm:1}.}

\begin{proof}[A proof of Main Theorem \ref{mthm:1}]

Let $k_0$ be an arbitrary integer satisfying $k_0>k$ as presented.  

Let us use $x:=(x_1,\cdots,x_k)$ for (local) coordinates for ${\mathbb{R}}^k$.

Let us use $y:=(y_1,\cdots,y_{k^{\prime}})$ for (local) coordinates for ${\mathbb{R}}^{k^{\prime}}$ where $k^{\prime}:=k_0-k$.

We take two steps to complete the proof. \\
\ \\
STEP 1 Defining a set in $M_D \subset {\mathbb{R}}^{k_0}$, which is a real algebraic hypersurface and a smooth regular compact submanifold of dimension $k_0-1$ with no boundary. \\

First we define $M_{D_0}:= \{(x,y) \in {\mathbb{R}}^k \times {\mathbb{R}}^{k^{\prime}}={\mathbb{R}}^{k_0} \mid  {\prod}_{j=1}^l (f_{{\rm P}, S_j}(x))-||y||^2=0\}$.\\
For the notation, remember the rule $||y||={\Sigma}_{j=1}^l {y_j}^2$.

We show this is also a smooth regular submanifold in ${\mathbb{R}}^{k_0}$. 
We consider the partial derivative of the function ${\prod}_{j=1}^l  (f_{{\rm P}, S_j}(x))-{\Sigma}_{j=1}^l {y_j}^2$ for variants $x_j$ and $y_j$. 
First we take a point $(x_0,y_0) \in M_{D_0}$ such that ${\prod}_{j=1}^l (f_{{\rm P}, S_j}(x_0))>0$. 
We use $x_0:=(x_{0,1},\cdots,x_{0,k})$ and $y_0:=(y_{0,1},\cdots,y_{0,k})$.
Here we consider the partial derivative of the function for some $y_j$ and we have the value $2{y_j}=2y_{0,j} \neq 0$. The differential of the restriction of the function ${\prod}_{j=1}^l  (f_{{\rm P}, S_j}(x))-{\Sigma}_{j=1}^l {y_j}^2$ at $(x_0,y_0) \in M_{D_0}$ is not of rank $0$ and this is not a singular point of the function.

Second we take a point $(x_{S_a},y_{S_a}) \in M_{D_0}$ such that $f_{{\rm P}, S_a}(x_{S_a})=0$. 
By the assumption on the hypersurfaces $S_b$ and the polynomials $f_{{\rm P}, S_j}(x_{S_b})$, $f_{{\rm P}, S_{a^{\prime}}}(x_{S_a})> 0$ for $a^{\prime} \neq a$. The polynomial function defined canonically from the polynomial $f_{{\rm P}, S_{a}}$ is assumed to have no singular points on $S_{a}$.
We use $x_{S_a}:=(x_{S_a,1},\cdots,x_{S_a,k})$ and $y_{S_a}:=(y_{S_a,1},\cdots,y_{S_a,k})$.
Here we consider the partial derivative of the function for some $x_{j}$ and we have the non-zero value
represented as the product of the partial derivative of the function $f_{{\rm P}, S_a}(x)$ for $x_{j}$ at $(x_{S_a},y_{S_a})$ and the product of $l-1$ numbers defined as the values of polynomials (or the canonically defined polynomial functions) in the family $\{f_{{\rm P}, S_j}\}_{j=1}^l$ at $x_{S_a}$ except the number $j \neq a$.  The differential of the restriction of the function at $(x_{S_a},y_{S_a}) \in M_{D_0}$ is not of rank $0$ and this is not a singular point of the function.

We have shown that $M_{D_0}$ is a smooth regular submanifold by the implicit function theorem. 

We define $M_{D}$ as the set of all points in $M_{D_0}$ such that $x \in \overline{D} \supset D$.
We investigate a small neighborhood of each point in $M_{D}$.

First we consider a point $p_1 \in D$ and a point $(p_1,q_1) \in M_{D}$ and take its sufficiently small open neighborhood $U_{p_1,q_1}$ in ${\mathbb{R}}^{k_0}$. For any point in $M_{D_0} \bigcap U_{p_1,q_1}$, by the definition, it is also a point in $M_{D}$.
Second we consider a point $p_2 \in \partial \overline{D}$ in the boundary $\partial \overline{D} \subset \overline{D}$ and a point $(p_2,q_2) \in M_{D}$ and take its sufficiently small open neighborhood $U_{p_2,q_2}$ in ${\mathbb{R}}^{k_0}$. Take an arbitrary
 point $(p^{\prime},q^{\prime})$ in $M_{D_0} \bigcap U_{p_2,q_2}$. By the definition and the assumption on the hypersurfaces $S_b$ and the polynomials $f_{{\rm P}, S_j}(x_{S_b})$, we can know $f_{{\rm P}, S_{b^{\prime}}}(p^{\prime})>0$ for $1 \leq b^{\prime} \leq l$ except for one $b^{\prime}:={b_0}^{\prime}$. $f_{{\rm P}, S_{{b_0}^{\prime}}}(p^{\prime})<0$ cannot occur by the form of the function ${\prod}_{j=1}^l  (f_{{\rm P}, S_j}(x))-{\Sigma}_{j=1}^l {y_j}^2$. $(p^{\prime},q^{\prime})$ is also a point in $M_{D}$.

We have shown that $M_{D}$ is also a smooth regular submanifold with $M_{D}=M_{D_0}$. By the form of the function ${\prod}_{j=1}^l (f_{{\rm P}, S_j}(x))-{\Sigma}_{j=1}^l {y_j}^2$ and the compactness of the closure $\overline{D}$, it is also a smooth compact manifold with no boundary. \\
\ \\
STEP 2 The composition of the restriction of the canonical projection mapping $(x_1,x_2) \in {\mathbb{R}}^{k_0}$ to $x_1 \in {\mathbb{R}}^k$ to the submanifold $M_D$ with the restriction of the given projection ${\pi}_{k,1}$ or ${\pi}_D$ in Definition \ref{def:1}. \\

First restrict the canonical projection mapping $(x_1,x_2) \in {\mathbb{R}}^{k_0}$ to $x_1 \in {\mathbb{R}}^k$ to the submanifold $M_D$. We thus have a surjection onto $\overline{D}$.
We restrict this to the preimage of $D$. This is, by the form of the
 function, regarded as a projection and a submersion. If we restrict this to the preimage of the boundary $\partial \overline{D}$, then, by the form of the
 function, we have a smooth and real algebraic embedding onto $\partial \overline{D}$. We compose the surjection onto $\overline{D}$ with the ${\pi}_D$ to have a new real algebraic function. \\
\ \\
By our definitions and situations, we can see that the composition obtained before is regarded as a desired function for $k_0>k+1$ where we need to respect connectedness of the preimages. 

We note some. Some may help us to understand our arguments more rigorously. First, our resulting function is a function for \cite{saeki2}, having finitely many singular values. Second, our map on $M_D$ onto $\overline{D}$ is, topologically, regarded as a so-called {\it special generic} map. The class of special generic maps contains Morse functions with exactly two singular points on spheres, or Morse functions in the so-called Reeb's theorem, and canonical projections of unit spheres. See \cite{saeki1} for fundamental theory on special generic maps and some advanced results on manifolds admitting such maps. For construction of special generic maps related to our construction of the map on $M_D$ onto $\overline{D}$, consult also the preprints \cite{kitazawa7, kitazawa11} of the author for example.
  
This completes the proof. 

.
	
\end{proof}

\begin{Ex}
	FIGURE 1 shows two simplest explicit cases. 
	
	The upper figure shows a Poincar\'e-Reeb graph for the pair of the algebraic domain surrounded by $l \geq 1$ circles centered at points and of fixed radii and a canonical projection into (a copy of) the $1$-dimensional real affine space where $l \geq 1$ is an arbitrary positive integer. It shows a graph with exactly $2$ vertices of degree $1$, exactly $2(l-1)$ vertices of degree $3$, and exactly $2(l-1)+l=3l-2$ edges.
	
	The lower figure shows a Poincar\'e-Reeb graph for the pair of the algebraic domain surrounded by $l \geq 1$ circles centered at points and of fixed radii and a canonical projection into (a copy of) the $1$-dimensional real affine space where $l \geq 2$ is an arbitrary integer greater than $1$. It shows a graph with exactly $2$ vertices of degree $1$, exactly $2$ vertices of degree $l+1$ and exactly $l+2$ edges.
\begin{center}
	\begin{figure}
		
		\includegraphics[height=45mm, width=40mm]{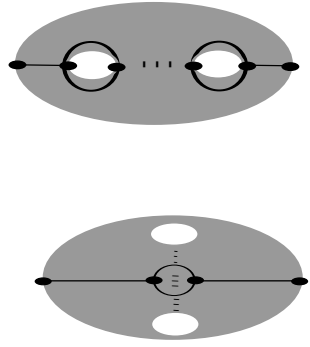}
		
		\caption{Some simplest Poincar\'e-Reeb graphs for Main Theorem \ref{mthm:1}. Small dots are for abbreviation of circles, edges and vertices.}
		\label{fig:1}

	\end{figure}
	\end{center}
	
\end{Ex}

We present remarks related to our result.
\begin{Rem}
	In the proof of Main Theorem \ref{mthm:1}, $||y||={\Sigma}_{j=1}^l {y_j}^2$ can be replaced by a polynomial of the form ${\Sigma}_{j=1}^l k_{1,j}{y_j}^{2k_{2,j}}$ with arbitrary positive integers $k_{1,j}$ and $k_{2,j}$ for example.
\end{Rem}

\begin{Rem}

According to the preprint \cite{bodinpopescupampusorea}, for any graph of some certain wide class, we can obtain some algebraic domain respecting the situation that the underlying $2$-dimensional real affine space and a more general projection are given. More precisely, we also have a Poincar\'e-Reeb graph for the pair of the real affine space and the general projection and the graph is isomorphic to the given graph as a graph. 
It tries to obtain domains arguing in the topological category or the class of $C^r$ which may not be the class of $C^{\infty}$ with $r \geq 1$ on the differentiability. After some arguments, it applies so-called Weierstrass-type theorem and approximations.

This can give various examples if the algebraic domains satisfy our conditions. However, investigating such conditions is in general difficult. See Example 2.2 and FIGURE 2 of the preprint as an explicit example.
\end{Rem}
We present another remark which is not directly related to our study in the present situation. We also wish this is closely related to our study in the future.  
\begin{Rem}
Here let ${\mathbb{C}}^k$ denote the $k$-dimensional complex space, whose underlying Euclidean space is of $2k$-dimensional and the real affine space of course. 
It is also a simplest complex algebraic manifold.
Let $\mathbb{C}:={\mathbb{C}}^1$.
It has been difficult to construct very explicit examples of real algebraic functions into $\mathbb{R}$ or maps into higher dimensional real affine spaces on explicitly given closed and connected real algebraic manifolds via explicit polynomial maps.
In \cite{sakurai}, Sakurai gives an explicit example via celebrating theory of Milnor on links of complex polynomials \cite{milnor}. He first considers the polynomial function on the $3$-dimensional complex space ${\mathbb{C}}^3$ mapping $(z_1,z_2,z_3) \in {\mathbb{C}}^3$ to ${z_1}^2+{z_2}^2+{z_3}^2 \in \mathbb{C}$ and the link associated with this link is represented as the intersection of the unit sphere $S^5$ in ${\mathbb{R}}^6$ and the zero set of the polynomial.  He restricts a very explicit complex linear function on the outer complex space ${\mathbb{C}}^3={\mathbb{R}}^6$ to the link, which is diffeomorphic to the $3$-dimensional real projective space, and obtains a smooth map into ${\mathbb{R}}^2$. This map enjoys nice properties. The image of the set of all singular points of the map is two smoothly and disjointly embedded circles. This is conjectured to be essentially same as a so-called {\it round} fold map in \cite{kitazawa1, kitazawa2, kitazawa3, kitazawa6, kitazawa10, kitazawasaeki1, kitazawasaeki2}, by Osamu Saeki and the author.
\end{Rem}
\section{Acknowledgement}
The author would like to thank Osamu Saeki and Shuntaro Sakurai for interesting discussions on \cite{sakurai} and construction of explicit smooth maps which are also real algebraic. This with an interesting talk in the conference "Singularity theory and its applications (RIMS-Sing 1)"\\
(http://www.math.kobe-u.ac.jp/HOME/saji/math/conf2022/spmon/index.html)\\
 by M. S. Sorea on \cite{bodinpopescupampusorea} and her interesting answers to questions on construction of nice domains and extensions of the presented study and results to graphs of wider classes by the author has motivated us to start studies on our present paper.


\begin{thebibliography}{25}
\bibitem{bodinpopescupampusorea} A. Bodin, P. Popescu-Pampu and M. S. Sorea, \textsl{Poincar\'e-Reeb graphs of real algebraic domains}, arXiv:2207.06871.
\bibitem{kitazawa1} N. Kitazawa, \textsl{On round fold maps} (in Japanese), RIMS Kokyuroku Bessatsu B38 (2013), 45--59.
\bibitem{kitazawa2} N. Kitazawa, \textsl{On manifolds admitting fold maps with singular value sets of concentric spheres}, Doctoral Dissertation, Tokyo Institute of Technology (2014).
\bibitem{kitazawa3} N. Kitazawa, \textsl{Fold maps with singular value sets of concentric spheres}, Hokkaido Mathematical Journal Vol.43, No.3 (2014), 327--359.
\bibitem{kitazawa4} N. Kitazawa, \textsl{On Reeb graphs induced from smooth functions on $3$-dimensional closed orientable manifolds with finitely many singular values}, Topol. Methods in Nonlinear Anal. Vol. 59 No. 2B, 897--912, arXiv:1902.08841.
\bibitem{kitazawa5} N. Kitazawa, \textsl{On Reeb graphs induced from smooth functions on closed or open surfaces}, Methods of Functional Analysis and Topology Vol. 28 No. 2 (2022), 127--143, arXiv:1908.04340.
\bibitem{kitazawa6} N. Kitazawa, \textsl{Round fold maps and the topologies and the differentiable structures of manifolds admitting explicit ones}, submitted to a refereed journal, arXiv:1304.0618.

\bibitem{kitazawa7} N. Kitazawa, \textsl{Notes on explicit special generic maps into Euclidean spaces whose dimensions are greater than $4$}, a revised version is submitted based on positive comments (major revision) by referees and editors after the first submission to a refereed journal, arxiv:2010.10078.
\bibitem{kitazawa8} N. Kitazawa, \textsl{On Reeb graphs induced from smooth functions on $3$-dimensional closed manifolds which may not be orientable}, a revised version is submitted to a refereed journal after based on positive comments by editors and referees after the second submission to a refreed journal, arXiv:2108.01300.
\bibitem{kitazawa9} N. Kitazawa, \textsl{Realization problems of graphs as Reeb graphs of Morse functions with prescribed preimages}, submitted to a refereed journal, arXiv:2108.06913.
\bibitem{kitazawa10} N. Kitazawa,\textsl{Round fold maps on $3$-dimensional manifolds and their integral and rational cohomology rings}, arXiv:2301.07008.
\bibitem{kitazawa11} N. Kitazawa, \textsl{A note on cohomological structures of special generic maps}, a revised version is submitted based on positive comments by referees and editors after the second submission to a refereed journal.
\bibitem{kitazawasaeki1} N. Kitazawa and O. Saeki, \textsl{Round fold maps on $3$-manifolds}, accepted for publication after a refereeing process and to appear in Algebraic \& Geometric Topology, arXiv:2105.00974.
			\bibitem{kitazawasaeki2} N. Kitazawa and O. Saeki, \textsl{Round fold maps of $n$-dimensional manifolds into ${\mathbb{R}}^{n-1}$}, submitted to a refereed journal, arXiv:2111.13510.
\bibitem{martinezalfaromezasarmientooliveira} J. Martinez-Alfaro, I. S. Meza-Sarmiento and R. Oliveira, \textsl{Topological  classification of simple Morse Bott functions on surfaces}, Contemp. Math. 675 (2016), 165--179.%
\bibitem{masumotosaeki} Y. Masumoto and O. Saeki, \textsl{A smooth function on a manifold with given Reeb graph}, Kyushu J. Math. 65 (2011), 75--84.
\bibitem{michalak} L. P. Michalak, \textsl{Realization of a graph as the Reeb graph of a Morse function on a manifold}. Topol. Methods in Nonlinear Anal. 52 (2) (2018), 749--762, arXiv:1805.06727.
\bibitem{milnor} J. Milnor, \textsl{Singular points of complex hypersurfacs}, Annals of Mathematics Studies, No. 61, Princeton University Press, Princeton, N. J.; University of Tokyo Press, Tokyo, 1968.
\bibitem{reeb} G. Reeb, \textsl{Sur les points singuliers d\'{}une forme de Pfaff compl\'{e}tement int\`{e}grable ou d\'{}une fonction num\'{e}rique}, Comptes Rendus
 Hebdomadaires des S\'{e}ances de I\'{}Acad\'{e}mie des Sciences 222 (1946), 847--849.
\bibitem{saeki1} O. Saeki, \textsl{Topology of special generic maps of manifolds into Euclidean spaces}, Topology Appl. 49 (1993), 265--293.
\bibitem{saeki2} O. Saeki, \textsl{Reeb spaces of smooth functions on manifolds}, International Mathematics Research Notices, maa301, https://doi.org/10.1093/imrn/maa301, arXiv:2006.01689.
\bibitem{sakurai} S. Sakurai, Master Thesis, Kyushu. Univ..
\bibitem{sharko} V. Sharko, \textsl{About Kronrod-Reeb graph of a function on a manifold}, Methods of Functional Analysis and
 Topology 12 (2006), 389--396.
\bibitem{sorea1} M. S. Sorea, \textsl{The shapes of level curves of real polynomials near strict local maxima},  Ph. D. Thesis, Universit\'e de Lille, Laboratoire Paul Painlev\'e, 2018.
\bibitem{sorea2} M. S. Sorea, \textsl{Measuring the local non-convexity of real algebraic curves}, J. Symbolic Compute. 109 (2022), 482--509.
\end{thebibliography}
\end{document}